\documentclass[11pt]{article}

\usepackage{graphicx}
\usepackage{amssymb}
\usepackage{mathtools} 
\usepackage{url}

\title{An implicit gradient-descent procedure for minimax problems}

\author{Montacer Essid \and Esteban G. Tabak  \and Giulio Trigila}

\begin{document}

\maketitle

\begin{abstract}
A game theory inspired methodology is proposed for finding a function's saddle points. 
While explicit descent methods are known to have severe convergence issues, implicit methods are natural in an adversarial setting, as they take the other player's optimal strategy into account. The implicit scheme proposed has an adaptive learning rate that makes it transition to Newton's method in the neighborhood of saddle points. Convergence is shown through local analysis and, in non convex-concave settings, thorough numerical examples in optimal transport and linear programming. An ad-hoc quasi Newton method is developed for high dimensional problems, for which the inversion of the Hessian of the objective function may entail a high computational cost. 
\end{abstract}

\section{Introduction}
Saddle point problems occur in a wide variety of applications ranging from mechanics and computational fluid dynamics to constrained optimization \cite{benzi2005numerical}, optimal transport \cite{ELT} and more recently to machine learning with the advent of generative adversarial neural network \cite{goodfellow2014generative}. The general structure of these problems is formulated 
in terms of the mini-maximization of a Lagrangian function:
\begin{equation}
  \min_x \max_y L(x,y), \quad x \in R^{n_x}, \quad y \in R^{n_y}.
  \label{minmax}
\end{equation}
The following is a list of examples directly related to the discussions below:

\medskip

\noindent
{\bf Example 1: equality-constrained minimization} 
  $$ \min_x f(x) \quad \hbox{subject to} \quad  g(x) = 0. $$
  Introducing Lagrange multipliers $y$ yields
  $$ \min_x \max_y L(x,y) = f(x) - y^t g(x). $$

Often some components of $x$ and $y$ are required to be non-negative:

\medskip

\noindent
{\bf Example 2: inequality-constrained minimization} 
  $$ \min_x f(x) \quad \hbox{subject to} \quad  g(x) \ge 0, $$
where introducing Lagrange multipliers $y$ yields
  $$ \min_x \max_{y\ge 0} L(x,y) = f(x) - y^t g(x). $$
  When there are both equality and inequality constraints, only the Lagrange multipliers $y_j$ attached to the inequalities are required to be non-negative.
  
\medskip

 \noindent 
{\bf Example 3: two-player zero-sum games} 
  
  $$ \min_x \max_y y^t A x, \quad \hbox{with } x,y \ge 0,  \quad \sum_i x_i = \sum_j y_j = 1. $$ 
 Introducing Lagrange multipliers $\lambda$ and $\mu$ for the equality constraints, yields
  $$ \min_{x \ge 0, \mu} \max_{y \ge 0, \lambda}\ L = {y}^t\ A\ x - \lambda \left(\sum_i x_i - 1\right) - \mu \left(\sum_j y_j  - 1\right). $$  

A more recent development formulates problems of interest as nonlinear adversarial games. Examples include generative adversarial networks \cite{goodfellow2014generative}, as well as the following: 

\medskip

\noindent
{\bf Example 4: adaptive optimal transport}
  
  $$ \min_{\alpha} \max_{\beta}\left[ \sum_i w^x_i\ g\left(\nabla\phi\left(x_i,\alpha\right), \beta\right) 
    - \sum_j w^y_j\ e^{g\left(y_j, \beta \right)}\right]. $$

%
%

%

%

\medskip

This article proposes a game-theory inspired methodology for the numerical solution of minimax problems, \emph{implicit twisted-gradient descent}. The adjective ``twisted'' refers to the fact that one player descends the gradient while the other ascends it. The adjective ``implicit'' specifies that the two players act simultaneously, descending (in $x$) and ascending (in $y$) the Lagrangian $L$ in an anticipatory manner, i.e. following the gradient of $L$ estimated at the values of $(x,y)$  resulting from the current step. For small learning rates $\eta$, each step of this procedure converges to regular (twisted) gradient descent, while for large $\eta$ it converges to a Newton step.


There are analogies between the work presented here and the proposal in \cite{mertikopoulos2018mirror}, which implements a twisted gradient descent by an ad hoc modification of the mirror descent method. This modification is based on predictable sequences in which at each time step a guess on the future direction of the gradient is made. 
In the methodology proposed here, the anticipation of the next gradient uses the Hessian, making it possible to leverage the extensive optimization literature on Newton's method. An example in this direction is the development of a quasi-Newton-like method presented in section \ref{sec:QN}. 

The plan of the article is as follows: after this introduction, section \ref{sec:TwistedGD} introduces the basic step of the procedure for a given learning rate $\eta$. Section \ref{sec:Conv} proves the procedure's local convergence. Section \ref{sec:learning} proposes an adaptive criterion for evolving the learning rate. Section \ref{sec:QN} develops a quasi-Newton-like methodology that bypasses the need to evaluate the Hessian of $L$ or to invert any matrix. Section \ref{sec:Constraints} extends the methodology to situation where some or all variables are required to be positive. Section \ref{sec:Examples} shows examples of numerical results. Finally, section \ref{concl} includes some concluding remarks.  

\section{Implicit twisted gradient descent}\label{sec:TwistedGD}

We consider first the case without positivity constraints:
$$
  \min_x \max_y L(x,y), \quad x \in R^{n_x}, \quad y \in R^{n_y}.
$$
A gradient descent step, ``twisted''  as required in a minimax setting, is given by
\begin{eqnarray*}
  x^{n+1} &=& x^n - \eta \nabla_x L\big|_{x^{n},y^{n}} \\
  y^{n+1} &=& y^n + \eta \nabla_y L\big|_{x^{n},y^{n}} ,
\end{eqnarray*}
where $\eta > 0$ is the learning rate: the players with strategy $x$ and $y$ seek to decrease and increase $L$ respectively, and do so following the direction of their components of the gradient of $L$.  For compactness, we introduce the following notation: 
\begin{equation}
  z =  \left(\begin{array}{r}
     x \\
     y
     \end{array}\right), \quad 
   G =
     \left(\begin{array}{r}
     \nabla_x L \\
     \nabla_y L
     \end{array}\right), \quad J = \left(\begin{array}{rr}
                   I_x & 0 \\
                   0 & -I_y \end{array} \right),
\end{equation}
where $I_x$ and $I_y$ are identity matrices of size $n_x$ and $n_y$ respectively.
Then the descent step reads
\begin{equation}
z^{n+1} = z^n  - \eta\ J \ G^n.    
\label{egd} 
\end{equation}
Yet such a procedure may fail to converge. Consider the simple example with $L=x y$,
which has the unique min-maximizer $x=y=0$. Here twisted gradient descent would yield
\begin{eqnarray*}
  x^{n+1} &=& x^n - \eta y^n \\
  y^{n+1} &=& y^n + \eta x^n,
\end{eqnarray*}
or
$$ z^{n+1} = \left(\begin{array}{rr}
                 1 & -\eta \cr
                 \eta & 1 \end{array} \right) z^n, $$
which diverges, since its eigenvalues $\lambda^{\pm} = 1 \pm i \eta$ have absolute value greater than $1$. Even in the limit of infinitesimally small values of $\eta$, the solution moves in circles around the origin, following the system of ODEs
\begin{eqnarray*}
  \dot{x} &=& \hskip-0.3cm -y \\
  \dot{y} &=& x.
\end{eqnarray*}

One could argue that, from a game-theory perspective, each player would not merely move following the local gradient, but would try to anticipate how the other player will move. This suggests a form of implicit twisted gradient descent:
\begin{equation}
z^{n+1} = z^n  -  \eta\ J\ G^{n+1}.     
\label{implicit}
\end{equation}
Notice that applying (\ref{implicit}) to the example above now yields
$$ \left(\begin{array}{rr}
                 1 & \eta \cr
                 -\eta & 1 \end{array} \right) z^{n+1} = z^n, $$
with unconditional convergence to $z=0$, at a convergence rate that grows unboundedly with the learning rate $\eta$.

Yet in general (\ref{implicit}) cannot be solved in closed form for $z^{n+1}$. Instead, we may approximate $G^{n+1}$ using
\begin{equation}
  G^{n+1} \approx G^n + H^n \left(z^{n+1} - z^n\right),
  \label{HqN}
\end{equation}
where $H$ is the Hessian
\begin{equation}
  H = 
     \left(\begin{array}{rr}
     L_{xx} & L_{xy} \\
     L_{yx} & L_{yy}
     \end{array}\right).   
\end{equation}
Under this approximation, the scheme in (\ref{implicit}) yields
%
$$
  z^{n+1}  = z^n  - \eta\ J\ \left(G^n + H^n \left(z^{n+1} - z^n\right)\right),
$$
%
which reduces to
%
%
\begin{equation}
z^{n+1} = z^n  - \eta \left(J + \eta H^n\right)^{-1} G^n.     
\label{igd}
\end{equation}
This is the basic updating step of the proposed algorithm.

Notice that, as $\eta \rightarrow \infty$, the update in (\ref{igd}) converges to the Newton step
%
$$
z^{n+1} = z^n  - \left(H^n\right)^{-1} G^n
$$
%
while, for small $\eta$, it approximates the explicit twisted gradient descent step in (\ref{egd}). 

This proposal leaves us with some tasks:
\begin{enumerate}
  \item Prove convergence of the algorithm,
  \item develop a scheme for updating the learning rate $\eta$ so as to accelerate convergence,
  \item develop a way to avoid inverting possibly large matrices and, if possible, avoid computing the Hessian altogether,
  \item extend the procedure to situation where some or all variables have positivity constraints, and
  \item show numerical examples of the algorithm at work.  
\end{enumerate}
These tasks are addressed in the following sections.

\section{Local convergence}\label{sec:Conv}

The minimax theorem \cite{du2013minimax} guarantees the existence of a mini-maximizer of $L$ under the assumptions that $L(x,y)$ is quasiconvex in $x$ for all $y$ and quasiconcave in $y$ for all $x$.  Under these conditions, 
$$   \min_x \max_y L(x,y) =  \max_y \min_x L(x,y). $$
For smooth functions $L(x,y)$ satisfying these assumptions, the mini-maximizer $z^*$ of $L$ is uniquely characterized by the first order conditions
$$ G\left(z^*\right) = 0. $$
More generally, any point $z^*$ where $G$ vanishes and $L(z^*)$ satisfies quasiconvexity/concavity locally, is a local mini-maximizer of $L$. If the quasiconvexity/concavity conditions are satisfied globally, or at least in a neighborhood of the local optimal $z^*$ that the algorithm does not leave, it is enough to show that $\|G\|$ decreases to zero in order to guarantee convergence of the algorithm. This is proven below.

Under the first-order approximation to the gradient in (\ref{HqN}), the procedure in (\ref{igd}) yields
\begin{equation}
  G^{n+1} = \left(I - \eta H^n \left(J + \eta H^n\right)^{-1}\right) G^n =  
  J \left(J + \eta H^n\right)^{-1} G^n,
  \label{update_G}
\end{equation}
guaranteeing convergence to $G = 0$ when
\begin{equation}
  {G^{n+1}}^t H^n J G^{n+1} = {\nabla_x L^{n+1}}^t L_{xx} \nabla_x L^{n+1} - 
  {\nabla_y L^{n+1}}^t L_{yy} \nabla_y L^{n+1} > 0,
  \label{GHG}
\end{equation}
since then
\begin{eqnarray}\label{eq:decrG}
 \left\| G^n \right\|^2 &=& \left\| G^{n+1} + \eta H^n J G^{n+1} \right\|^2 \nonumber \\
 &=& \left\| G^{n+1} \right\|^2 + 
2\eta {G^{n+1}}^t H^n J G^{n+1} + \eta^2 \left\|H^n J G^{n+1}\right\|^2 \nonumber \\
 &>&  \left\| G^{n+1} \right\|^2.
 \label{Gn1}
\end{eqnarray}
Notice that the condition in (\ref{GHG}) is weaker than the hypothesis of convexity in $x$ and concavity in $y$, as these are only required to hold in the direction of the gradient and, moreover, one of the two may not hold provided that the other compensates in the sum. Hence, for problems that do not satisfy the convexity requirement globally, one must add additional checks to exclude convergence to non-optimal points. For example, the Lagrangian
$$ L(x,y) = \frac{1}{2} \left(x^2 + y^2\right) $$
has no finite min-maximizer, yet the algorithm will converge to $(0,0)$ if started at $(a,0)$ for any value of $a$. On the other hand, this is an unstable trajectory: when the starting point is moved to $(a, \epsilon)$ with $\epsilon \ne 0$, the solution necessarily diverges for $\eta$ smaller than $2$, as we show below.

Notably, for $\eta$ large enough, the inequality (\ref{Gn1}) is satisfied under a condition much weaker than (\ref{GHG}), namely that $H^n J G^{n+1} \ne 0$, i.e. that $G^{n+1} \ne G^n$. In particular, if one adopts a quasi-Newton method, where (\ref{HqN}) is satisfied by construction, and let $\eta \rightarrow \infty$, we have guaranteed convergence to $G = 0$ for any $L(x,y)$. Thus, under certain conditions, it may be convenient to adopt a learning rate $\eta$ as large as possible. This will guide our choice for selecting a learning rate below.

\section{Determination of the learning rate}\label{sec:learning}

In order to turn the implicit gradient descent (\ref{igd}) into an algorithm, one needs a mechanism to decide at each step which learning rate $\eta$ to use.

From the arguments above, once close enough to the optimum, one should increase $\eta$ as much as possible so as to accelerate convergence, with $\eta = \infty$ yielding Newton's method. Yet Newton's method is blind to whether one is minimizing or maximizing the objective function. In our minimax context, it could converge to points where $G=0$ that are not minima over the $x$ and maxima over the $y$. 
Therefore, a mechanism to reject too large values of $\eta$ is required. Unlike in pure minimization scenarios, we cannot use the decrease of the objective function as an acceptance test. However, a simple extension applies: every step should satisfy 
\begin{equation}
 L\left(x^{n+1},y^{n}\right) \le L\left(x^{n+1},y^{n+1}\right) \le L\left(x^{n},y^{n+1}\right).
 \label{Check}
\end{equation}
This agrees with the anticipatory game idea underlying the method: given $y^{n+1}$, the player with strategy $x$ should make sure to decrease $L$, and given $x^{n+1}$, the player with strategy $y$ should make sure to increase $L$. Thus a step not satisfying (\ref{Check}) should be rejected.



In subsection \ref{3ex} we quantify the effects of this constraint on three prototypical examples for which we can write a closed form for (\ref{igd}) and (\ref{Check}).

\medskip

\subsection{A strategy for evolving the learning rate}

Rather than updating the learning rate $\eta$ directly, we update a surrogate $\mu$, and then build $\eta$ dividing $\mu$ by $\|G^n\|$, so as to normalize the step-size in $z$-space. The algorithm proposed is the following:

\begin{enumerate}

\item Set an initial guess $z^0$ and an initial value $\mu^0$. 

\item At each step, update $\mu$ through $\mu^{n+1}=\min\left(\alpha \mu^n, \mu_{max}\right)$, with $\alpha > 1$, $\mu_{max} \gg 1$. Update $z^n$ to $z^{n+1}$ through (\ref{igd}) with $\eta=\mu/\|G^n\|$. If the conditions in (\ref{Check}) are
not satisfied, reduce $\mu^{n+1}$ (for instance halving it) until either they are satisfied or $\mu^{n+1}$ is smaller than a prescribed threshold.

\item Stop when either $\|G^{n+1}\|$ is smaller than a prescribed threshold or the number of steps reaches a prescribed maximum.

\end{enumerate}

\subsection{Three examples}
\label{3ex}


In this subsection we consider three simple prototypical examples where a closed expression for the constraints in (\ref{Check}) can be derived:

\begin{enumerate}
  \item $L = xy$,
  
  \item $L = \frac{x^2 - y^2}{2}$,
  
  \item $L = \frac{x^2 + y^2}{2}$.

\end{enumerate}

The first represents a saddle point not satisfying the convexity conditions $L_{xx} > 0$, $L_{yy} < 0$, yet having a global solution ($x=y=0$), the second does satisfy these conditions globally, and the third has no solution, so we would like $y$ to blow up: with no local minimax solution, the algorithm should explore other areas of $(x,y)$-space. 

\subsubsection{$L = xy$}

Gradient, Hessian, $J$:

$$ G = \left(\begin{array}{c}
           y \\
           x \end{array} \right), \quad
     H = \left(\begin{array}{cc} 
     0 & 1 \\
     1 & 0 \end{array} \right), \quad
     J = \left(\begin{array}{rr} 
     1 & 0 \\
     0 & -1 \end{array} \right).
           $$
The update rule (\ref{igd}) then yields
$$   \left(\begin{array}{c}
           x \\
           y \end{array} \right)^+ = \left(\begin{array}{c}
           x \\
           y \end{array} \right) - \eta \left(J + \eta H\right)^{-1} \left(\begin{array}{c}
           y \\
           x \end{array} \right)  = \frac{1}{1+\eta^2} \left(\begin{array}{c}
           x - \eta y \\
           y + \eta x \end{array} \right),
 $$
where $(x,y)$ stands for the current (i.e. $n$th) state, and $(x,y)^+$ for the next one.
Notice that here the larger $\eta$ the better, as increasing $\eta$ brings us closer to the solution $(0,0)$. We have

$$L(x^+,y) = \frac{(x-\eta y) y}{1+\eta^2}, \quad L(x,y^+) = \frac{(y+\eta x) x}{1+\eta^2}, $$

$$ L(x^+,y^+) = \frac{(y+\eta x) (x-\eta y)}{(1+\eta^2)^2}.$$

Notice that
$$ L(x^+,y^+) - L(x^+,y) = \frac{\eta}{(1+\eta^2)^2} (x-\eta y)^2 $$
and
$$ L(x,y^+) - L(x^+,y^+) = \frac{\eta}{(1+\eta^2)^2} (y+\eta x)^2,$$
both non-negative for all positive $\eta$, hence imposing no restrictions. This is in line with the fact that, in this case, the solution of the problem can be reached in just one step by adopting $\eta = \infty$.


\subsubsection{$L = \frac{1}{2} \left(x^2 - y^2\right)$}

$$ G = \left(\begin{array}{r}
           x \\
           -y \end{array} \right), \quad
     H = \left(\begin{array}{rr} 
     1 & 0 \\
     0 & -1 \end{array} \right), \quad
     J = \left(\begin{array}{rr} 
     1 & 0 \\
     0 & -1 \end{array} \right).
           $$
$$   \left(\begin{array}{c}
           x \\
           y \end{array} \right)^+ = \left(\begin{array}{c}
           x \\
           y \end{array} \right) - \eta \left(J + \eta H\right)^{-1} \left(\begin{array}{r}
           x \\
           -y \end{array} \right)  = \frac{1}{1+\eta} \left(\begin{array}{c}
           x \\
           y \end{array} \right).$$
Again, the larger $\eta$ the better. We have

$$L(x^+,y) = \frac{1}{2} \left(\frac{x^2}{(1+\eta)^2}-y^2\right), \quad
L(x,y^+) = \frac{1}{2} \left(x^2-\frac{y^2}{(1+\eta)^2}\right), $$

$$ L(x^+,y^+) = \frac{1}{2} \frac{1}{(1+\eta)^2} \left(x^2-y^2\right). $$

Then both
$$ L(x^+,y^+) - L(x^+,y) = \frac{1}{2} \frac{(1+\eta)^2-1}{(1+\eta)^2} y^2 $$
and
$$ L(x,y^+) - L(x^+,y^+) = \frac{1}{2} \frac{(1+\eta)^2-1}{(1+\eta)^2} x^2 $$
are automatically non-negative for positive $\eta$, thus imposing no constraints. Once again this is in line with the fact that the exact solution can be reached in one step by adopting $\eta = \infty$.


\subsubsection{$L = \frac{1}{2} \left(x^2 + y^2\right)$}

Here the update rule (\ref{igd}) yields
\begin{equation}
  \begin{pmatrix}
    x \\
    y
   \end{pmatrix}^{+}=
 \begin{pmatrix}
 \frac{1}{1+\eta} & 0 \\ 
 0 & \frac{1}{1-\eta}
 \end{pmatrix}
 \begin{pmatrix}
 x\\ 
 y
 \end{pmatrix} ,
\end{equation}
which, for $\eta > 2$, converges to $(x,y)=(0,0)$, which does not min-maximize $L$, despite having zero gradient. On the other hand,
the conditions in (\ref{Check}) are
$$ L(x^+,y^+) - L(x^+,y) = \frac{y^2}{2} \left(\frac{1}{(1-\eta)^2} - 1 \right) \ge 0,$$
$$ L(x^+,y^+) - L(x,y^+) = \frac{x^2}{2} \left(\frac{1}{(1+\eta)^2} - 1 \right) \le 0.$$
While the second of these imposes no constraint, the first restricts 
$\eta$ to be smaller than $2$, thus guaranteeing divergence. This is the required output for this problem with no minimax. In a more general setting, this divergence would correspond to leaving the region where $L$ does not have the right convexity/concavity in $x$ and $y$, hence opening the search for true minimax solutions elsewhere.


\section{Quasi implicit twisted gradient descent}\label{sec:QN}

The leading computational costs of the proposed procedure are the computation of the Hessian (which may not even be available in closed form) and the inversion of its mollified version, i.e. the calculation of the matrix
$$ B = \left(J + \eta H\right)^{-1}. $$
In the spirit of quasi-Newton methods \cite{nocedal1980updating}, one can replace $B$ with an estimation that is updated at each time-step using our knowledge of the gradient at two consecutive times, $G^n$ and $G^{n+1}$, since (\ref{update_G}) reads:
$$
  J G^{n+1} = B^n G^n.
$$
Of course, at the time of updating $B^n$, one does not yet know $G^{n+1}$. Instead, one can update $B$ correcting $B^{n-1}$ into a $B^*$ that would have satisfied this constraint at the prior step:
$$
  J G^{n} = B^* G^{n-1}, \quad B^* \rightarrow B^n.
$$

A significant difference with regular quasi-Newton methods though is that $B$ is not positive definite, unlike the Hessian in minimization problems. Thus, even though one could propose the equivalent to the BFGS recipe:
\begin{equation}
  B^{*} =  W_n^t B^{n-1} W_n + \frac{J G^{n} (J G^{n})^t}{(G^{n-1})^t J G^{n}},
  \label{BFGS_B}
\end{equation}
where
$$ W_n = I - \frac{G^{n-1}(J G^{n})^t}{(G^{n-1})^t J G^{n}}, $$
this could yield an uncontrollable large correction, since the denominator can vanish even for an arbitrarily small learning rate $\eta$, for which $B=J$ and $G^{n}=G^{n-1}$.

An alternative is to perform the rank-one update
%
$$  B^{n} = B^* = B^{n-1} +  \frac{\left(J G^{n}-B^{n-1} G^{n-1}\right)\left(J G^{n}-B^{n-1} G^{n-1}\right)^t}{(G^{n-1})^t \left(J G^{n}-B^{n-1} G^{n-1}\right)}. $$
Similarly to (\ref{BFGS_B}), this corrects $B$ so that it gives the right answer on $G^{n-1}$.
In order to avoid singularities when the denominator vanishes, we may write this as
\begin{equation}
  B^{*} = B^{n-1} + \alpha \frac{\left(J G^{n}-B^{n-1} G^{n-1}\right)\left(J G^{n}-B^{n-1} G^{n-1}\right)^t}{\left\|(J G^{n}-B^{n-1} G^{n-1}\right\|^2},
  \label{r2B}
\end{equation} 
with
$$ \alpha = \frac{\left\|(J G^{n}-B^{n-1} G^{n-1}\right\|^2}{(G^{n-1})^t \left(J G^{n}-B^{n-1} G^{n-1}\right)} $$
replaced with
$$ \alpha^* = \hbox{sign}(\alpha) \min(|\alpha|, \|B^{n-1}\|),$$
i.e. limiting the norm of the rank-one update to that of $B^{n-1}$. Applying a similar solution to eliminate possible singularities to (\ref{BFGS_B}) would have been problematic, as we would have had to fix not only the second term of the sum on the right hand side of (\ref{BFGS_B}) but also the matrix $W_{n}$.

\section{Inequality constraints}\label{sec:Constraints}

Often some or all $z_i$ are required to be in some subset, typically to be non-negative. We can limit consideration to this latter case with little loss of generality, since any constraint of the form $g(z) \ge 0$ can be reduced to the positivity of the corresponding Lagrange multiplier. So we have the problem
$$ \min_x \max_y L(x,y), \quad z(I) \ge 0, $$
where $z=\{x,y\}$, and $I$ indexes the subset of variables required to be non-negative. There are a number of ways to extend the procedure of this article to the case with inequalities; be discuss below two alternative methodologies:

\subsection{Change of variables}

The simplest way to enforce positivity without altering the algorithm is to make a change of variables that ensures positivity, for instance setting 
%
$$ L^*(x,y) = L(X(x), Y(y)), $$
where 
\begin{equation}
  X(x) = \begin{cases}
                        x & \hbox{for unrestricted variables} \cr
                         x^2 & \hbox{for variables required to be non-negative}
              \end{cases}
              \label{x_to_X}
 \end{equation}
 and similarly for $Y(y)$. This yields the unconstrained minimax problem
 $$ \min_x \max_y L^*(x,y) $$
 to which the procedure can be applied, and whose solution, once transformed into $Z=\{X,Y\}$, solves the original problem
 $$ \min_X \max_Y L(X,Y), \quad Z(I) \ge 0. $$

 A word of caution is in order though: the fact that, for $i \in I$, $z_i = 0 \Rightarrow L^*_{z_i} = 0$ creates potential suboptimal points where the procedure might stop. For instance, in constrained optimization problems, the Lagrange multipliers corresponding to inactive constraints are zero at the solution, but one often encounters along the way to the true solution, domains where some constraints that will be active in the final solution are temporarily inactive. Hence these Lagrange multipliers $z_i$ may reach machine zero values, at which point the corresponding derivatives of $L^*$ vanish. Because of this, these $z_i$ may fail to leave zero when the corresponding constraints become active again.
 
This issue can be addressed through a simple procedural change: after every step, compute the gradient of the original Lagrangian, i.e. $L_Z$, for the variables $\{Z_i\}$ that are close to zero, i.e. $\|z_i\| \le \epsilon$. Since $z_i$ should detach from zero  when the original gradient $L_{Z_i}$ pushes $Z_i={z_i}^2$ to be positive, we compute
$$ R_i = \max(-J_i^i L_{Z_i}, 0) $$
and update $z$ via
\begin{equation}
  z_i \rightarrow \tilde{z_i} = \sqrt{{z_i}^2 + \eta_0 R_i},
\end{equation}
where $\eta_0$ is a suitably small additional learning rate, restricted so as to satisfy the requirements in (\ref{Check}) between the states $z_i$ and $\tilde{z_i}$. To do this, we start with an arbitrary value for $\eta_0$ and reduce it, for instance by halving, until (\ref{Check}) is satisfied.
 
\subsection{Evolving barriers}

 A more conventional approach to handling positivity constraints is to add a logarithmic barrier:
\begin{equation}
 L(x,y) \rightarrow L^t(x,y) = L(x,y) + \frac{1}{t} \left[\sum_j \log(y_j) - \sum_i \log(x_i) \right].
\end{equation}
Here we can either solve the problem for an increasing sequence of values of $t$, adopting as initial values of $(x,y)$ for each subproblem their terminal values from the prior one, or take this to the limit, evolving $t$ smoothly at each step of the algorithm. 

%
%

\section{Examples}
\label{sec:Examples}

We illustrate the procedure through three examples: a simple two-dimensional one designed to illustrate the effects of non-convexity on quasi-implicit descent and the need for the constraints imposed on the learning rate, a linear programming one to illustrate the handling of inequality constraints when very many are simultaneously active, and an optimal transport problem to show a nonlinear adversarial example of current interest.  

\subsection{A two-dimensional example}
\begin{figure}[ht!]
  \begin{center}
      \begin{tabular}{ll}                                                                                      
      \resizebox{59mm}{!}{\includegraphics{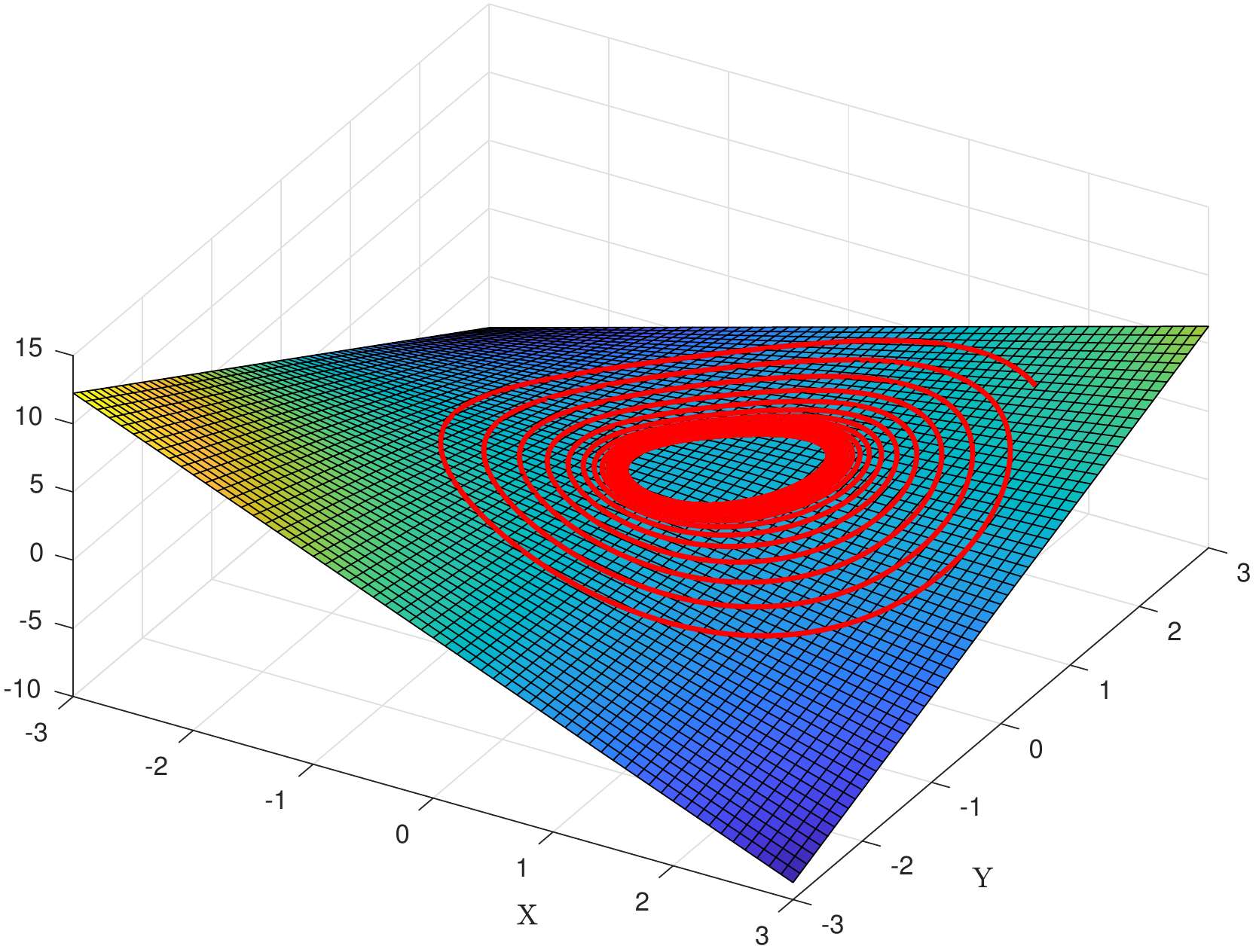}}&
      \resizebox{59mm}{!}{\includegraphics{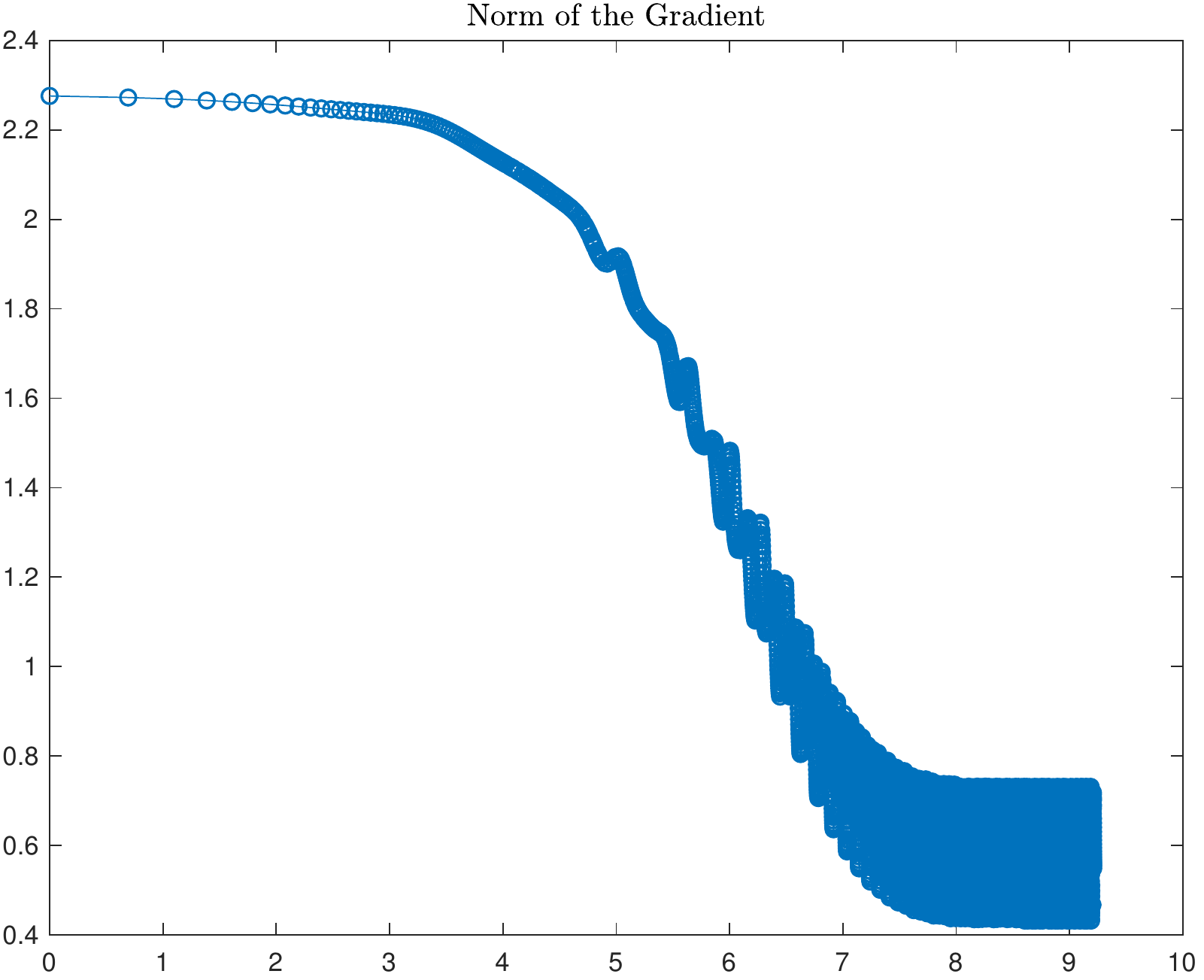}}\\
 		     \resizebox{59mm}{!}{\includegraphics{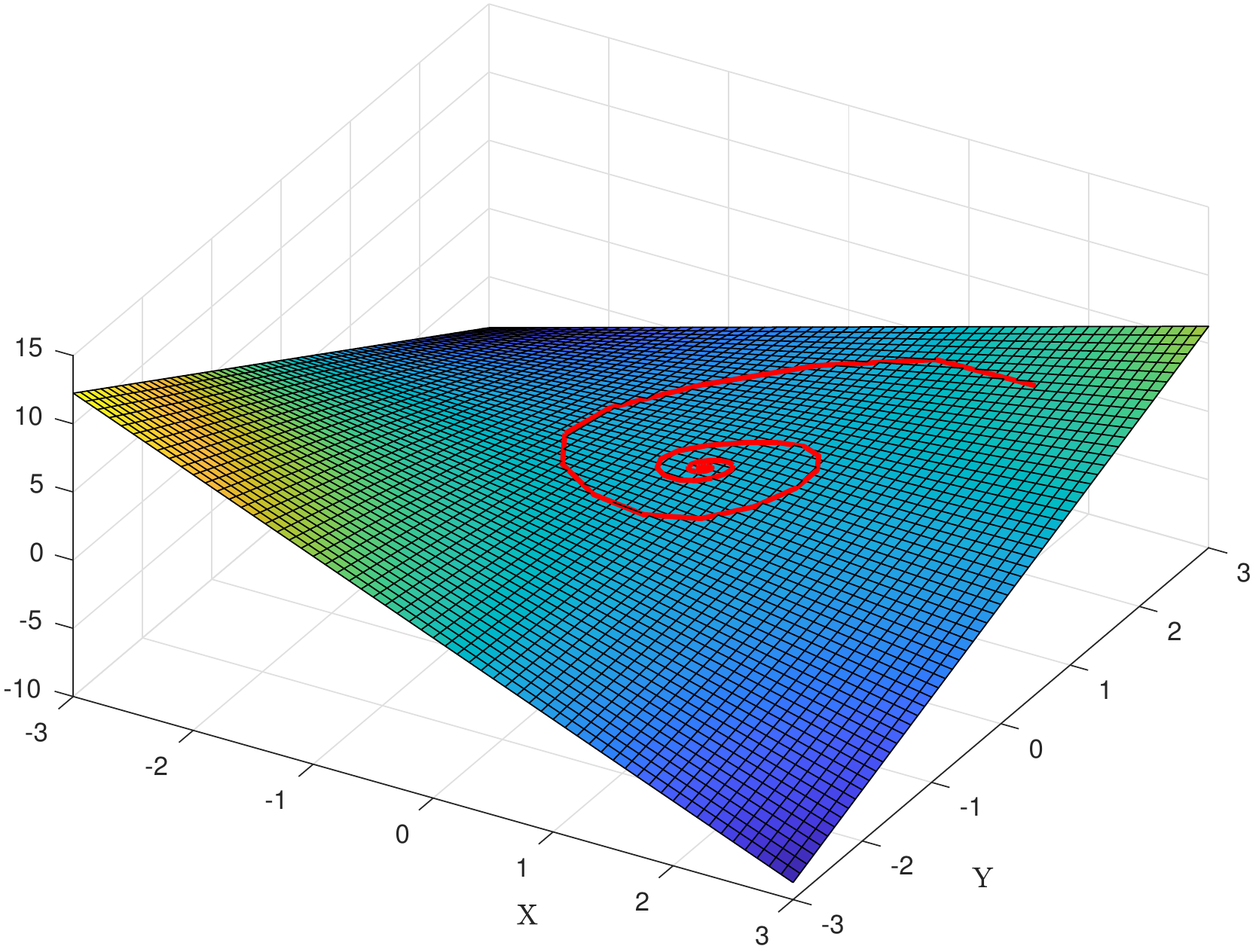}}&
      \resizebox{59mm}{!}{\includegraphics{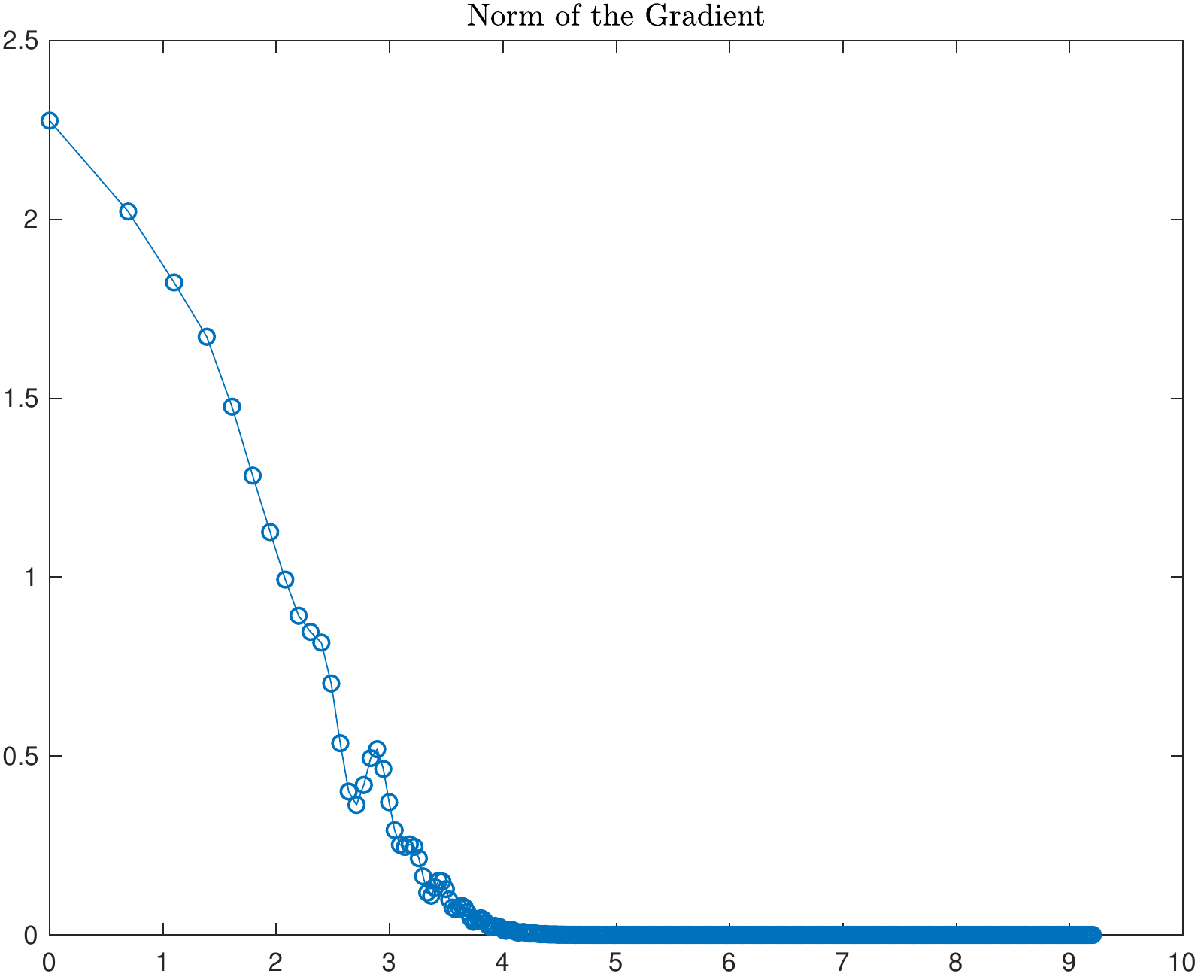}}\\
      \resizebox{59mm}{!}{\includegraphics{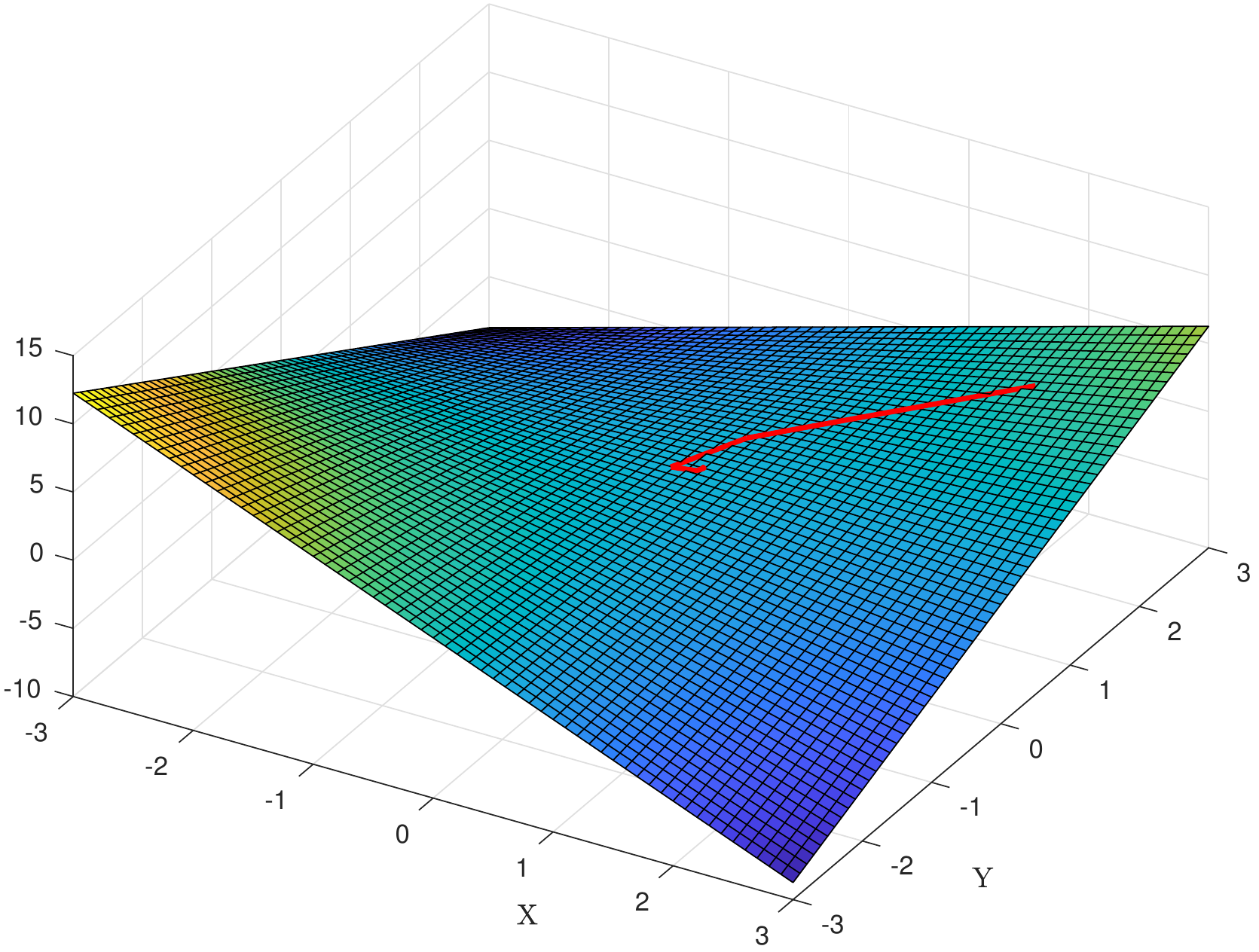}}&
      \resizebox{59mm}{!}{\includegraphics{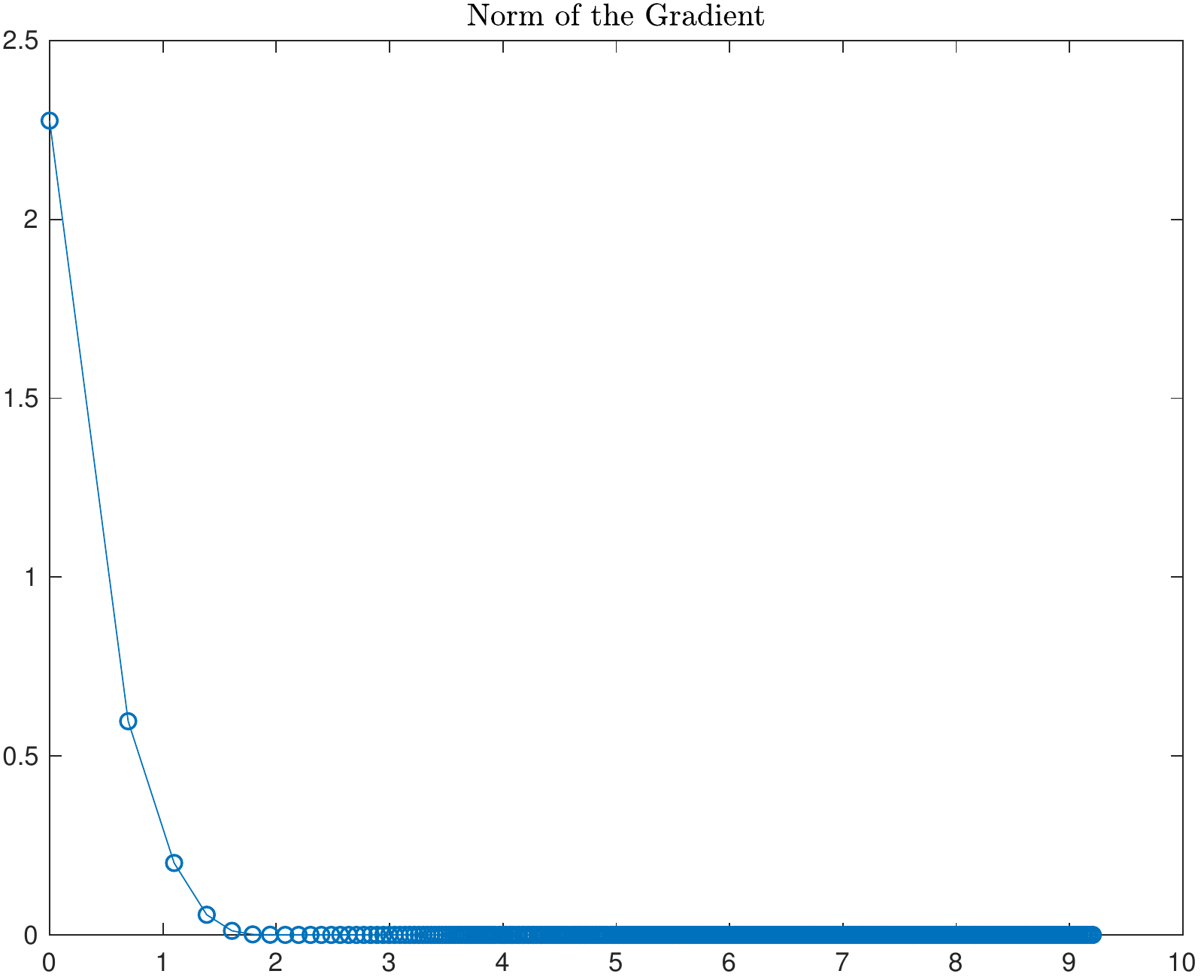}}\\
     \end{tabular}
       \end{center}
    \caption{Three trajectories using the algorithm in (\ref{igd}) with 3 different values of fixed $\eta$ to compute the saddle point of (\ref{simplesaddle}). The left column shows the trajectory and the right column the value of the norm of the gradient appearing in (\ref{igd}) as a function of the logarithm of the iteration step. Each row of the plot is obtained with values of $\eta$ equal to 0.05, 0.5 and 5 respectively. For a too small value of the ``looking forward" time $\eta$ the algorithm behaves essentially as the analogous gradient ascent-descent resulting in a periodic orbit. As the value of $\eta$ increases the gradient decreases as described by (\ref{eq:decrG}).}
    \label{fig:simple}
\end{figure}

\
This sub-section displays a numerical example of the implicit gradient descent and the quasi Newton method on a non-monotone saddle point problem (i.e. one in which the objective function is not convex-concave in the variables in which we are minimizing and maximizing respectively). The Lagrangian is
\begin{equation}\label{simplesaddle}
L(x,y)=(x-0.5)(y-0.5)+\frac{1}{3}e^{(-(x-0.5)^2 -(y-0.75)^2 )}.
\end{equation}
This function has a saddle point near $(0.5, 0.5)$ and a local maximum near $(0.5, 0.75)$. 
It has been observed in \cite{mertikopoulos2018mirror} that, in this case, first-order descent methods result in periodic orbits. Figure \ref{fig:simple} shows that this is indeed the case if the look ahead time $\eta$ in (\ref{igd}) is very small, effectively reducing (\ref{igd}) to an explicit algorithm. It is also interesting to notice that, in line with the discussion in section \ref{sec:Conv}, for large values of $\eta$ we reach a very fast convergence since, close to the saddle point of $f(x,y)$, the implicit (\ref{igd}) is essentially exact.
 
\begin{figure}[ht!]
  \begin{center}
      \begin{tabular}{ll}                                                                                      
      \resizebox{59mm}{!}{\includegraphics{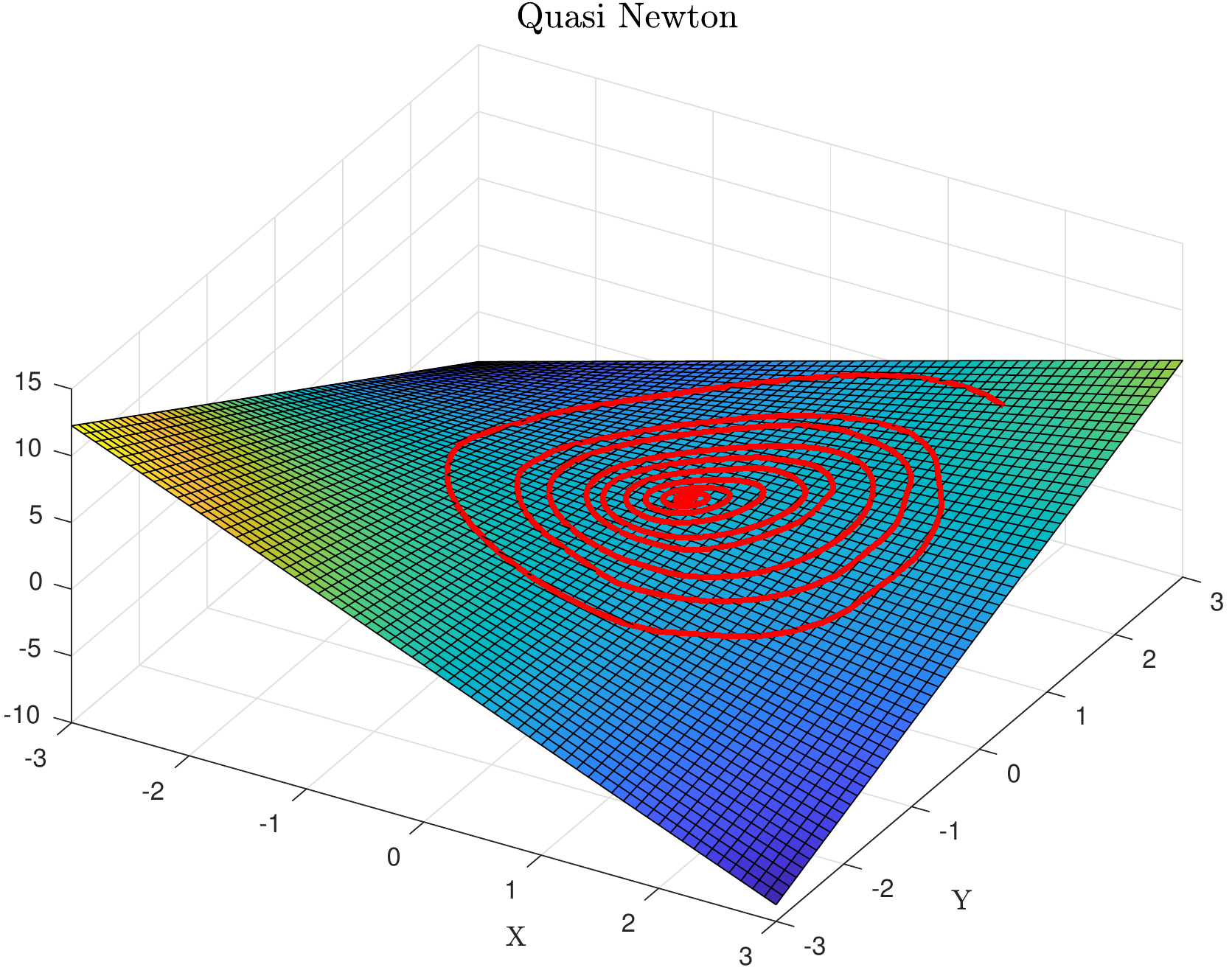}}&
      \resizebox{59mm}{!}{\includegraphics{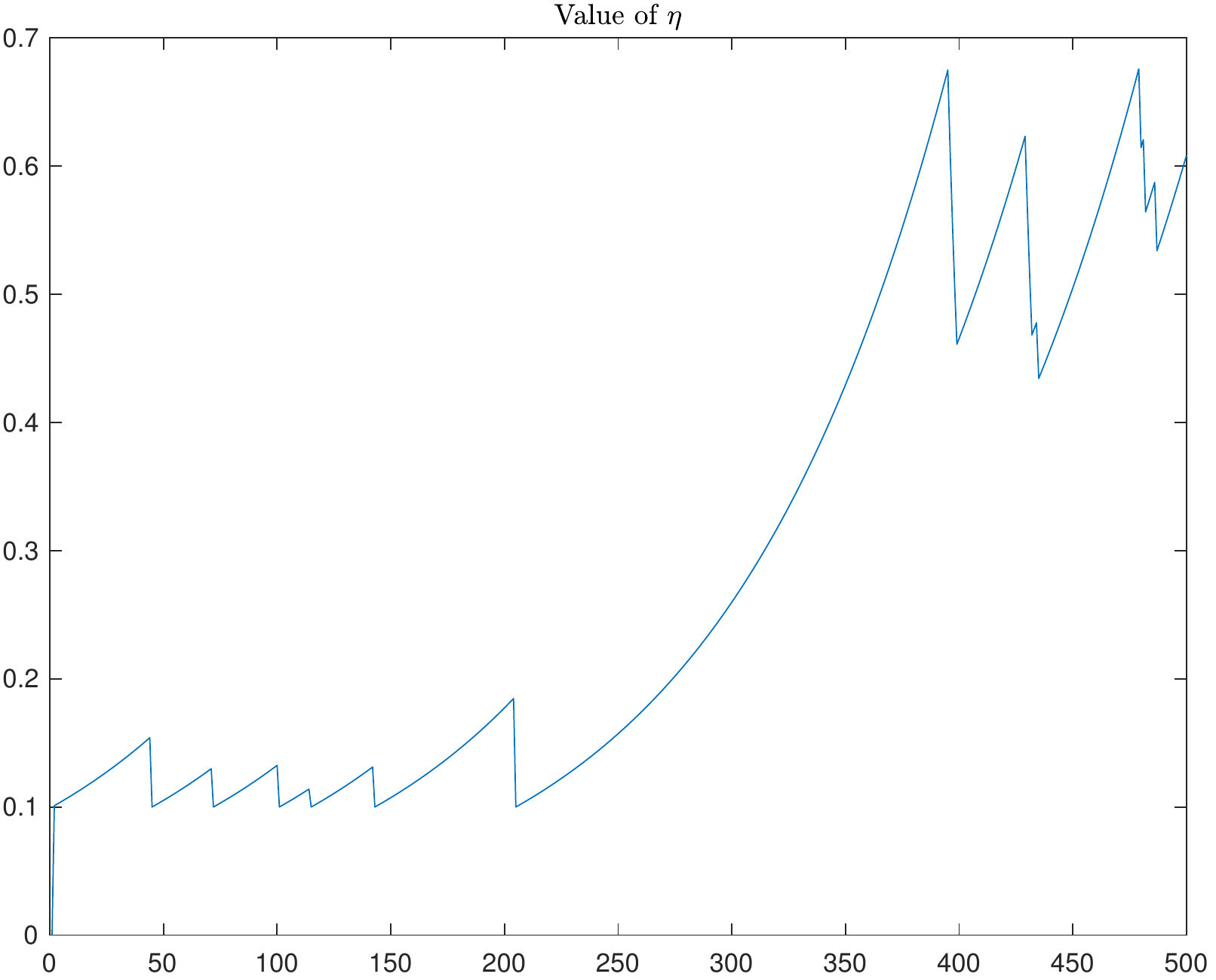}}\\
     \end{tabular}
    \caption{Trajectory obtained when using the Quasi Newton algorithm with variable $\eta$ as described in section \ref{sec:QN}. It can be seen that the learning rate $\eta$ get smaller in certain points of the trajectory due to the enforcing of the conditions in (\ref{Check}).}
    \label{fig:qn}
  \end{center}
\end{figure}

Figure \ref{fig:qn} shows the performance of the Quasi Newton algorithm with variable learning rates $\eta$ of section (\ref{sec:QN}). We see that the algorithm effectively ``learns'' the Hessian, leading to convergence. The jumps of the value of $\eta$ correspond to violation of the constraint in  (\ref{Check}).

\subsection{Linear programming}

We consider the standard linear programming problem
\begin{equation}
  \min_{X \ge 0} c^t X, \quad
  AX \ge b,
  \label{LP}
\end{equation}
which, introducing Lagrange multipliers $Y$ for the constraints, adopts the Lagrangian form
\begin{equation}
  \min_{X \ge 0} \max_{Y \ge 0}\ L^*(X,Y) = c^t X - Y^t \left(AX - b\right).
  \label{L*}
\end{equation}
To eliminate the positivity constraints, we introduce unconstrained variables $x$ and $y$ through $X = x^2$, $Y = y^2$, both understood component-wise, which yields the unconstrained minimax problem
\begin{equation}
  \min_{x} \max_{y}\ L(x,y) = c^t X(x) - Y(y)^t \left(AX(x) - b\right).
  \label{L}
\end{equation}

We have
\begin{equation}
   L_x = 2\left(c - A^t Y\right) .* x , 
   \quad L_y = 2\left(b - A X\right) .* y ,
\end{equation}
and
\begin{eqnarray*}     
L_{xx} = 2\ \hbox{diag}\left(c - A^t Y\right) , &&  L_{xy} = -4\ \hbox{diag}(x) A^t \hbox{diag}(y) \\
L_{yx} = -4\ \hbox{diag}(y) A\ \hbox{diag}(x) && L_{yy} = 2\ \hbox{diag}\left(b - A X\right).
\end{eqnarray*}
where the symbol `$.*$' denotes component-wise multiplication, and `$\hbox{diag}(x)$' denotes a diagonal matrix with the vector $x$ on its diagonal.

 \begin{figure}[ht!]
  \begin{center}
      \begin{tabular}{lll}                
      \hskip-1.4cm \resizebox{60mm}{!}{\includegraphics{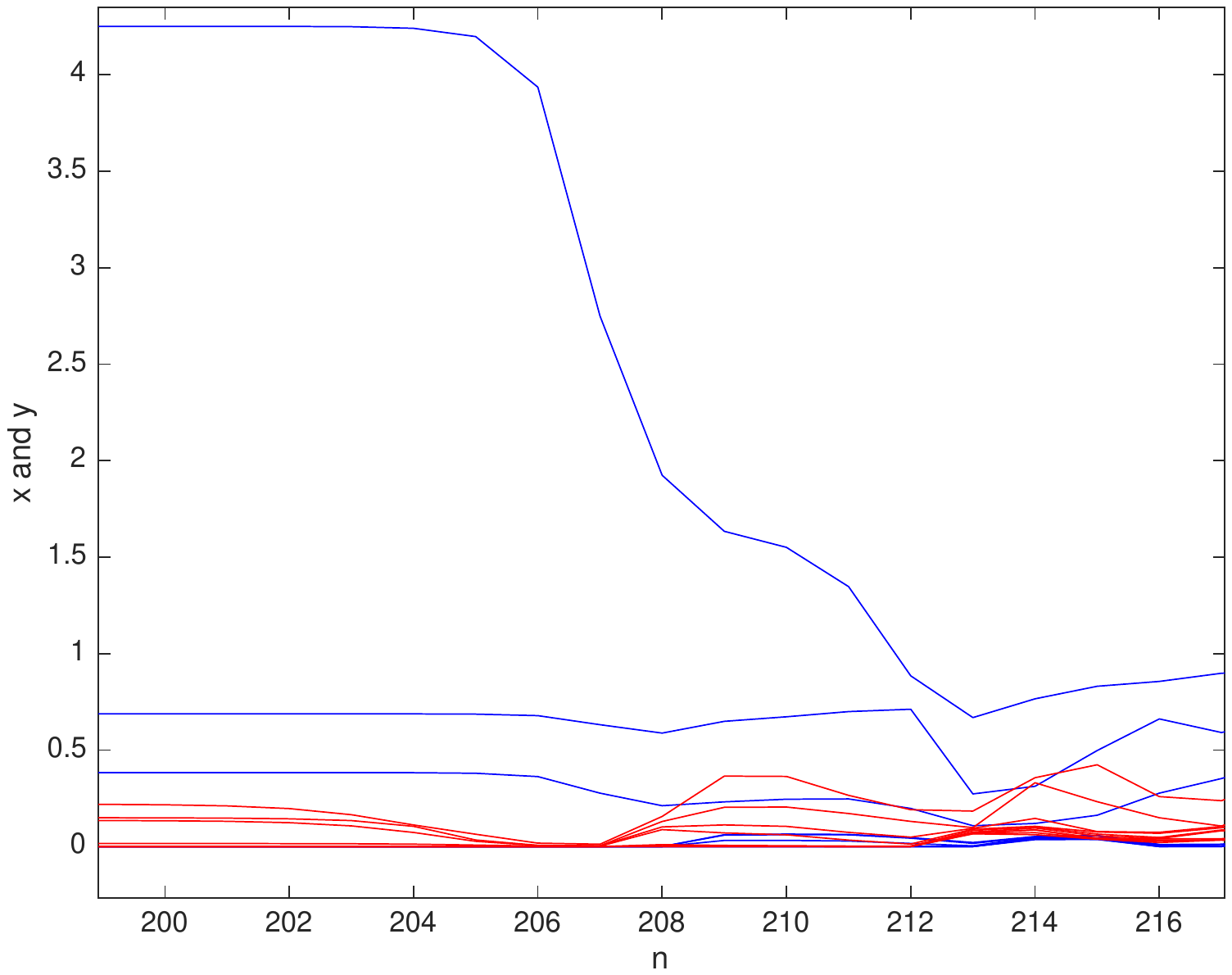}} &                                                                      
      \hskip-1.4cm\resizebox{60mm}{!}{\includegraphics{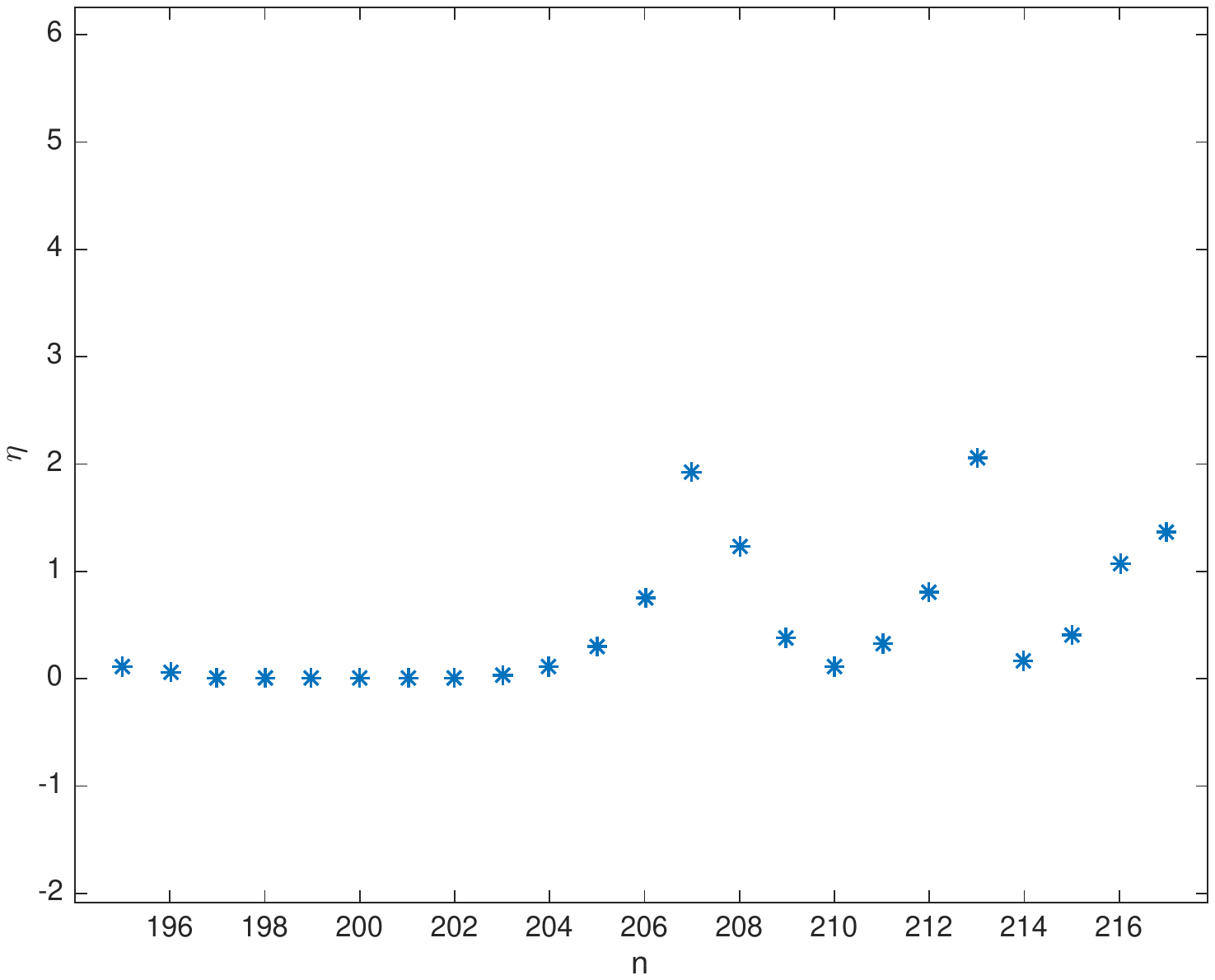}} &
     \hskip-1.4cm \resizebox{60mm}{!}{\includegraphics{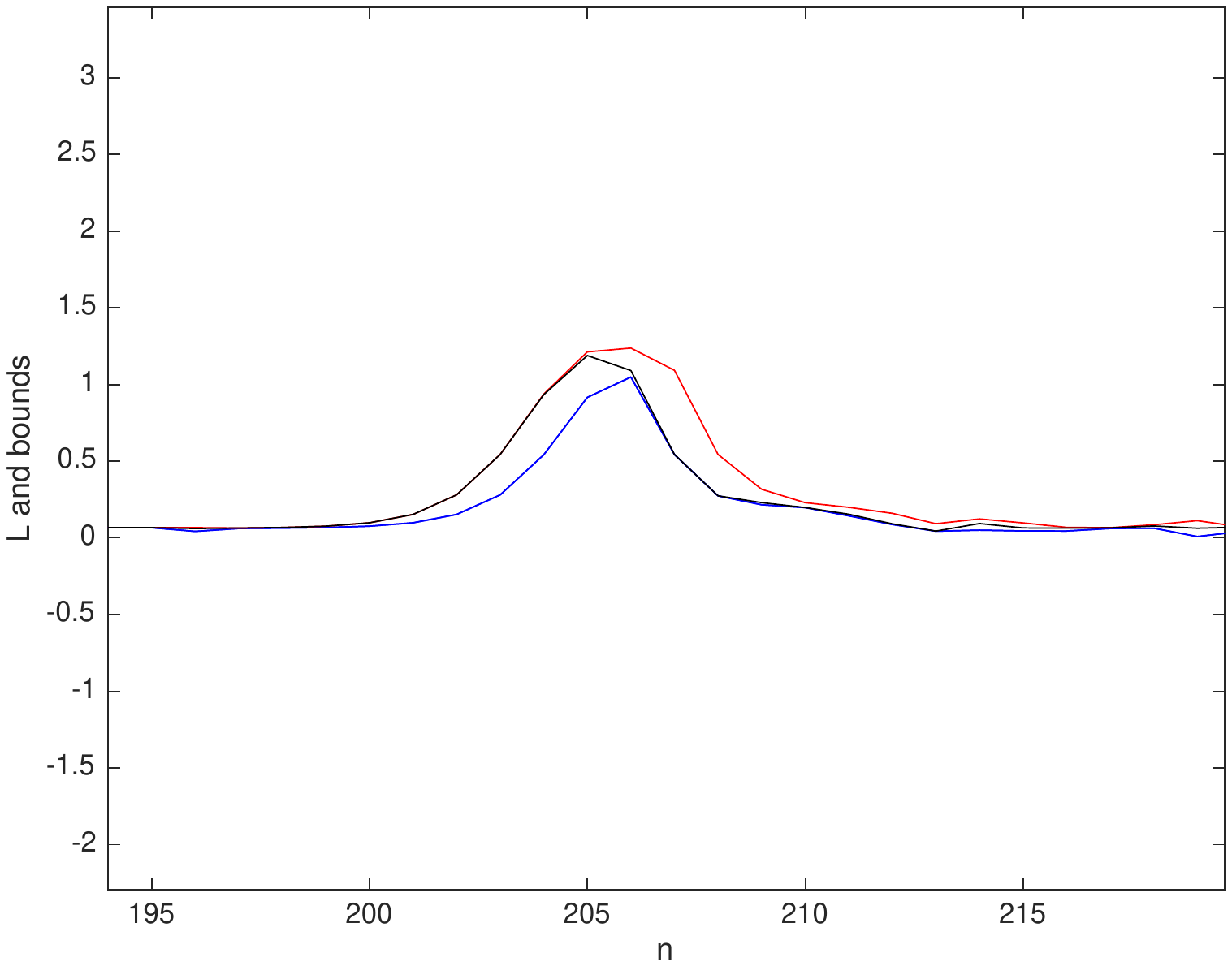}}
     \end{tabular}
    \caption{Linear Programming, zoom of the evolution near times ($n=206, 212$) where the active set changes considerably. Left panel: evolution of the 10 largest $x_i$ (in blue) and $y_j$ (in red). Middle panel: learning rate $\eta$. Right panel: the three values of $L$ appearing in the checks in (\ref{Check}), with the actual future $L$ in black and its upper and lower bounds in red and blue respectively.}
    \label{fig:LP}
  \end{center}
\end{figure}

For the example displayed in figure \ref{fig:LP}, we chose $n_x=117$, $n_y=114$. All entries of the matrix $A$ and the vectors $b$ and $c$ were drawn independently from the uniform distribution in $[0,1]$, thus guaranteeing feasibility.

The only free parameters of the procedure are the maximum learning rates, which we fixed at $10^7$, the rate $\alpha = 5.1$ at which $\mu$ is updated, and the initialization of $x$ and $y$, for which we picked quite arbitrarily
$$ x_0(1:n_x) = \sqrt{\frac{0.8}{n_x}}, \quad y_0(1:n_y) = \sqrt{\frac{0.4}{n_y}}. $$
For every realization of the problem, the procedure converges invariably to the right answer in 200-300 steps. A characteristic of this problem is that most positivity constraints are active, not only in the final solution but also at intermediate steps. Figure \ref{fig:LP} displays the working of the procedure at times where the active set changes significantly. We can see a local increase of the learning rate, corresponding to the opening of a significant gap between the lower and upper bounds for $L$ in (\ref{Check}).

\subsection{Optimal Transport}
 An adaptive, adversarial methodology has been developed in \cite{ELT} for the optimal transport problem \cite{monge1781memoire,kantorovich1942v}, between two distributions $\mu$ and $\nu$, known only through a finite set of independent samples. The problem consists in finding a global map $T$, pushing samples generated by the source $\mu$, so that their final distribution matches $\nu$, the one underlying the samples of the target. In addition, this map should minimize a transportation cost. For quadratic cost functions, the map $T$ must be given by the gradient $\nabla \phi$ of a convex potential $\phi$. We generate $T$ by composing many elementary non-linear functions $u_k$. Each of these $u_k$ minimizes a local optimal transport problem between two nearby samples $(x_i^{(k)})_{i=1,\ldots,n}$ and $(y_j^{(k)})_{j=1,\ldots,m}$. A global iterative procedure using displacement interpolation guarantees convergence to the unique optimizer.

In order to find these local non-linear maps $u$, we minimize the Kullback-Leibler divergence between the distributions underlying $u(x_i)$ and $y_j$. A variational characterization of the Kullback-Leibler divergence gives rise to the following formulation of the local problem:
\begin{equation}
\min_{u=\nabla\phi} \max_g \left\lbrace \frac{1}{n} \sum_i g(u(x_i)) - \frac{1}{m} \sum_j e^{g(y_j)} \right\rbrace
\end{equation}
The above mini-maximization can be interpreted as a two player game between the map $u$ and the lens $g$: as $u$ does its best to push the $x_i$'s toward the $y_j$'s, $g$ will focus on the areas where the mass transport has not yet been well achieved. This forces $u$ to correct those areas, and $g$ to find new locations requiring more work.

The maps $u$ and $g$ are parameterized using finite dimensional vectors $\alpha$ and $\beta$, and the problem is reduced to:
\begin{equation}
\min_{\alpha} \max_{\beta} \left\lbrace \frac{1}{n} \sum_i g_{\beta}(u_{\alpha}(x_i)) - \frac{1}{m} \sum_j e^{g_{\beta}(y_j)} \right\rbrace \equiv \min_{\alpha} \max_{\beta} L(\alpha, \beta)
\end{equation}

We solve each of those local optimal problems using the methodology described in this manuscript.

\begin{figure}[ht!]
  \begin{center}
      \begin{tabular}{ll}                                                                                      
      \resizebox{59mm}{!}{\includegraphics{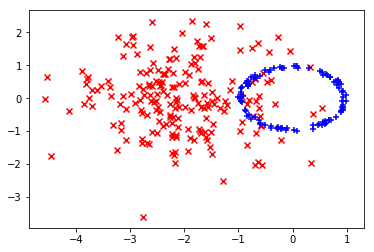}}&
      \resizebox{59mm}{!}{\includegraphics{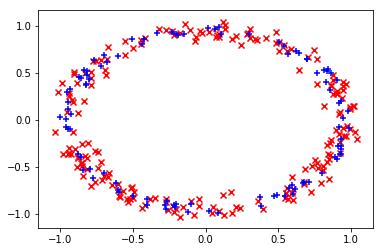}}\\
     \end{tabular}
    \caption{Initial and final configuration of the global optimal transport algorithm. The blue crosses represent the samples $(y_j)$, the red crosses on the left figure represent the samples $(x_i)$, and the red crosses in the right figure represent the samples generated by $T(x_i)$ where $T$ is a solution of the optimal transport algorithm}
    \label{fig:g2a}
  \end{center}
\end{figure}

\begin{figure}[ht!]
  \begin{center}
      \begin{tabular}{ll}                                                                                      
      \resizebox{120mm}{!}{\includegraphics{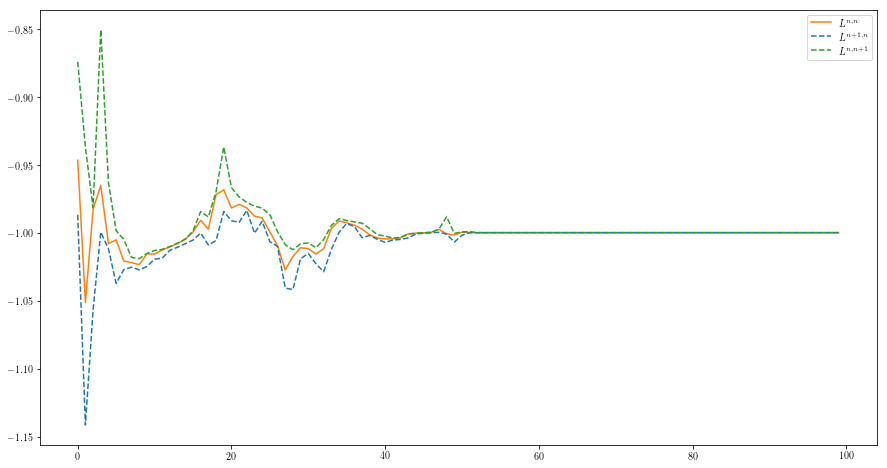}}
     \end{tabular}
    \caption{Values of the Lagrangian at $L(\alpha^n, \beta^n) \equiv L^{n,n}, L(\alpha^{n+1},\beta^n) \equiv L{n+1,n}$ and $L(\alpha^n, \beta^{n+1}) \equiv L^{n,n+1}$ for the last local optimal transport problem of the first global iteration.}
    \label{fig:optimalTransportL}
  \end{center}
\end{figure}

Figure \ref{fig:g2a} presents the original configuration of samples and the result of the global procedure, applied to data $\{x_i\}$ drawn from a Gaussian and $\{y_j$\} from the uniform distribution on the perimeter of a circle.
Figure \ref{fig:optimalTransportL} displays the objective function at each step, for the last local optimal transport problem of the first global iteration. In addition to $L(\alpha^n, \beta^n)$, displayed in orange, the upper bound $L(\alpha^n, \beta^{n+1})$ and the lower bound $L(\alpha^{n+1}, \beta^n)$ are displayed in green and blue respectively.

\section{Conclusions}
\label{concl}
This article presents an implicit twisted gradient descent strategy for the numerical computation of saddle points. Explicit methods are by nature non-anticipatory, which makes them often fail to converge, ending out in periodic or outward spiraling orbits around a saddle point. Instead, the algorithm proposed here is implicit, or anticipatory from a game theory perspective, as each player includes their adversary's best strategy in their own planning. This is proved to yield local convergence, which acquires a super-quadratic rate as the learning rate grows and the methodology converges to Newton's. The strategy proposed for updating the learning rate is consistent with the anticipatory nature of the algorithm: the rate should grow rapidly near saddle points, but is bounded by the requirement that, given the adversary's choice, each player should be improving their game. This guarantees convergence near saddle points and local divergence near stationary points of the Lagrangian that do not solve the minimax problem, points toward which regular Newton would otherwise converge.

The use of an implicit algorithm requires the inversion of a matrix, which can be quite large for high-dimensional problems. To alleviate this computational cost, the analogue of a quasi-Newton formulation of the algorithm is developed, which updates directly the inverse $B$ of the mollified Hessian at the core of the algorithm. This not only serves the purpose of avoiding matrix inversion, but also eliminates the need to compute or estimate second derivatives of the Lagrangian. 

Numerical tests are performed on three representative problems: a small-dimensional minimax problem that does not satisfy global convex-concavity, linear programming with a high number of inequality constraints, and a recently proposed adversarial methodology for optimal transport. In their diversity, they illustrate the versatility of the proposed methodology, which can be applied without modifications to virtually any minimax problem. It has been the author's experience that having such a general tool at one's disposal encourages the formulation of problems of interest in adversarial terms, a natural characterization that one would otherwise often avoid for lack of a straightforward methodology for their numerical solution.

\section{Acknowledgements}
The work of E. G. Tabak  was partially supported by NSF grant DMS-1715753 and ONR grant N00014-15-1-2355. The work of M. Essid was partially supported by NSF grant DMS-1311833.

\appendix
\addcontentsline{toc}{chapter}{APPENDICES}

\section{Appendix: Further corrections on the quasi Newton step} 

If required, one may refine the updating step by considering also the Hessian $H$. Since (\ref{HqN}) reads
$$
  \Delta G = H^n \Delta z, \quad \hbox{where} \quad \Delta G = G^{n+1} - G^n \quad
  \hbox{and} \quad \Delta  z = \left(z^{n+1} - z^n\right),
$$
the Hessian can also be updated via a rank-one update:
\begin{equation}
  H^{n+1} = H^n + \beta^* \frac{\left(\Delta G-H^n \Delta z\right)\left(\Delta G-H^n \Delta z\right)^t}{\|\left(\Delta G-H^n \Delta z\right)\|^2},
  \label{r1H}
\end{equation} 
$$ \beta^* = \hbox{sign}(\beta) \min(|\beta|, \|H^n\|),$$
$$ \beta = \frac{\|\left(\Delta G-H^n \Delta z\right)\|^2}{\left(\Delta z\right)^t\left(\Delta G-H^n \Delta z\right)}.$$
Hence we can independently update $B$ and $H$. Yet $B$ and $H$ are not independent, as they satisfy
\begin{equation}
 (J + \eta H) B = I.
 \label{Identity}
\end{equation}
We can use (\ref{Identity}) to improve our current estimations of $H$ and $B$, performing one step of gradient descent for
\begin{equation}
 F(B,H)=\frac{1}{2}\|(J + \eta H) B - I\|^2.
 \label{Identity_OF}
\end{equation}
Introducing a vector $v$ with all the entries of $H$ and $B$, we have $F=F(v)$, and we can write
$$ v = v - \nu \nabla_v F, $$
with
$$ \nu = \frac{2 F}{\|\nabla_v F\|^2}.$$

For reference, here goes the explicit calculation of the gradient:

$$  A = (J + \eta H) B - I, $$

$$ \frac{\partial A_m^n}{\partial H_i^j} = \eta \delta_i^m B_j^n, \quad  
\frac{\partial F}{\partial H_i^j} = \eta \sum_n A_i^n B_j^n, \quad \frac{\partial F}{\partial H} = \eta A B^t,$$ 

$$ \frac{\partial A_m^n}{\partial B_i^j} = \delta_j^n \left(J +\eta H\right)_m^i, \quad  
\frac{\partial F}{\partial B_i^j} = \sum_n A_m^j \left(J +\eta H\right)_m^i, \quad \frac{\partial F}{\partial B} = \left(J +\eta H\right)^t A.$$ 

The most expensive component of this procedure is the computation of the norms required to estimate $\nu$. To simplify this, notice that
$$ \nu = \frac{2 F}{\|\nabla_v F\|^2}= \frac{\|A\|^2}{\|\eta A B^t + \left(J +\eta H\right)^t A\|^2} \le \frac{1}{1 + \eta^2 \left(\|B\|^2+\|H\|^2 \right)}.$$
The norms of $B$ and $H$ can be estimated via
$$ \|B\|^2 \approx \|B u\|^2, \quad \|H\|^2 \approx \|H v\|^2, $$
where
$$ u^{n+1} = \frac{B u^n}{\|B u^n\|}, \quad u^{n+1} = \frac{H v^n}{\|H v^n\|}.$$

If at any point one would like to attempt a Newton step (i.e. set $\eta=\infty$), one needs access to $H^{-1}$. This can be updated in an entirely similar fashion to $B$ and $H$, since
$$ H^{-1} \Delta G = \Delta z, $$
which can be enforced through the rank-one update
\begin{equation}
  \left(H^{-1}\right)^{n+1} = \left(H^{-1}\right)^n + \frac{\left[\Delta z - \left(H^{-1}\right)^n \Delta G \right](\Delta G)^t}{(\Delta G)^t \Delta G}.
  \label{r1H3}
\end{equation} 
Again, the consistency between $H$ and $H^{-1}$ can be enforced through the gradient descent of
$$ \frac{1}{2}\|H^{-1} H - I\|^2. $$

Descending a matrix $A$ toward the inverse of another matrix $E$ involve the product of matrices, and operation that, while not inexpensive, can be easily parallelized.


%
%

\bibliographystyle{plain}      
\bibliography{ImplMM}
   
\end{document}